\documentclass[12pt]{amsart}
\usepackage{amssymb}
\usepackage{verbatim}
\usepackage{mathrsfs}
\usepackage{mathtools}

\usepackage{xcolor}

\usepackage{caption}

\usepackage{blindtext}

\newtheorem{thm}{Theorem}[section]

\numberwithin{equation}{section}

\providecommand{\arcosh}{\operatorname{arcosh}}

\newcommand{\alt}{\alpha}
\newcommand{\Lo}{ l_1 }

\newcommand{\HyG}{ {}_2F_1 }

\newcommand{\Mod}[1]{\ (\textup{mod}\ #1)}

\mathtoolsset{showonlyrefs}

\begin{document}

\title[Asymptotics of a Gauss hypergeometric function ]{Asymptotics of a Gauss hypergeometric function related to moments of symmetric-square $L$-functions I}

\author[O.Balkanova]{Olga Balkanova}
\address{
Steklov Mathematical Institute of Russian Academy of Sciences, 8 Gubkina st., Moscow, 119991, Russia}
\email{balkanova@mi-ras.ru}

\begin{abstract}
We prove an asymptotic formula for a special case of the Gauss hypergeometric function, which arises in the explicit formula for the second moment of Maass form symmetric-square $L$-functions.
\end{abstract}


\maketitle


\section{Introduction}
This paper is devoted to proving an asymptotic formula for the function
\begin{equation}\label{2f1(1/4 3/4 1)}
\HyG\left(1/4+ir(1-\alt),3/4+ir(1-\alt),1+2ir;z\right),\quad 0<z<1,\,
r\to\infty,\,\alpha=t_0/r\to0.
\end{equation}
This is one of the two special cases of the Gauss hypergeometric function arising in the explicit formula for the first moment of Maass form symmetric square $L$-functions on the critical line $\rho=1/2+2it_0$ (see \cite[Theorem 5, Lemma 7]{Bal} for details). The second special case is treated separately in \cite{Frol}. These two results, combined with the explicit formula, allow us to prove a uniform asymptotic formula for the considered moment, as shown in \cite{BalFrolEdin}.

\begin{thm}\label{prop: 2F1 asympt1}
Assume that $0<z<1,$ $r\to\infty$ and $r^{-1+\delta}<\alt< r^{-\delta}$. Then
\begin{multline}\label{2F1 main asymptI}
\HyG\left(1/4+ir(1-\alt),3/4+ir(1-\alt),1+2ir;z\right)=\\=
\frac{\exp\left(2ir\Lo(\alt,z)\right)}{(1-(1-\alt^2)z)^{1/4}}
\left(1+\sum_{j=1}^{N}\frac{c_j(\alt,r,z)}{(\alt r)^j}+O\left(\frac{\sqrt{r}}{(r\alt)^{N}}\right)\right),
\end{multline}
where $c_j(\alt,r,z)\ll1$ and
\begin{equation}\label{2F1 main asymptI Lo func}
\Lo(\alt,z)=\log2-\alt\log(1+\alt)-\log(1+\sqrt{1-(1-\alt^2)z})+\alt\log(\alt+\sqrt{1-(1-\alt^2)z}).
\end{equation}
\end{thm}
This asymptotic formula generalizes \cite[Lemma 4.6]{BF2mom},  where the case $\alt=0$ has been considered. Note that for $\alt=0$ and $z=4/x^2$ the main term in \eqref{2F1 main asymptI}
reduces to the expression (see \cite[4.37.19]{HMF})
\begin{equation}
\frac{\exp\left(2ir\left(\log2-\log(1+\sqrt{1-4/x^2})\right)\right)}{(1-4/x^2)^{1/4}}=
x^{2ir}\exp(-2ir\arcosh{x/2})
\left(\frac{x^2}{x^2-4}\right)^{1/4},
\end{equation}
which coincides with the main term in \cite[Lemma 4.6]{BF2mom}.
\section{The first case: $0<z<1-\delta$}
We start the proof of \eqref{2F1 main asymptI} by considering the case $0<z<1-\delta$.  In these settings, several methods may be used to derive the required asymptotic formula, including the Liouville-Green method  (by generalizing \cite[Lemma 2.4]{Zav} and \cite[Corollary 7.3]{BF2}) or a standard  version of the saddle point method (generalizing \cite[Lemma 20]{Bal}). However, the common problem of these two methods is that they does not seem to be applicable when $z$ is very close to $1$. Let us provide some explanation related to this issue starting with the Liouville-Green method and later discussing the saddle point method.

The Liouville-Green method is based on studying solutions of differential equations, see \cite[Ch. 6, pp. 10--12]{Olver}.   More precisely, using \cite[p.96, f.(7)-(9)]{BE} we conclude that
the function
\begin{equation}\label{Y def}
Y(z)=z^{1/2+ir}(1-z)^{1/2-ir\alt}\HyG\left(1/4+ir(1-\alt),3/4+ir(1-\alt),1+2ir;z\right)
\end{equation}
satisfies the differential equation
\begin{equation}\label{Y difeq}
Y''(z)=\left((ir)^2h(\alt,z)+g(\alt,z)\right)Y(z),
\end{equation}
where
\begin{equation}
h(\alt,z)=\frac{1-(1-\alt^2)z}{z^2(1-z)^2},\quad
g(\alt,z)=-\frac{1-3z/4+3z^2/4}{4z^2(1-z)^2}.
\end{equation}
The  function $h(\alt,z)$ has double poles at the points $0$ and $1$ and a turning point at $z_0=(1-\alt^2)^{-1}$. If $z$ is not to close to $1$, then  this turning point does not play any role, and consequently, it is possible to apply the theory developed in  \cite[Ch.6 Sec.5.3, Ch.10 Sec. 4.1]{Olver}.  According to \cite[Ch.10, Eqs. (1.01), (2.01)-(2.03), (2.09), (2.10)]{Olver}
and \cite[Ch.10, Thms. 3.1 and 4.1]{Olver} one has
\begin{equation}\label{Y to W1W2}
Y(z)h^{1/4}(\alt,z)=C_1W_{n,1}(r,\xi)+C_2W_{n,2}(r,\xi),
\end{equation}
where
\begin{equation}\label{xi def}
\xi=\int\sqrt{h(\alt,z)}dz,
\end{equation}
and
\begin{equation}\label{W1W2 def}
W_{n,1}(r,\xi)=e^{ir\xi}\sum_{j=0}^{n-1}\frac{A_{j}(z)}{(ir)^j}+O(r^{-n}),\quad
W_{n,2}(r,\xi)=e^{-ir\xi}\sum_{j=0}^{n-1}\frac{(-1)^jA_{j}(z)}{(ir)^j}+O(r^{-n}),
\end{equation}
\begin{equation}\label{Aj def}
A_{j+1}(z)=\frac{-A'_j(z)}{2\sqrt{h(\alt,z)}}+\int\Lambda(z)A_j(z)dz,
\end{equation}
\begin{equation}\label{Lambda def}
\Lambda(z)=\frac{16g(\alt,z)h^2(\alt,z)+4h(\alt,z)h''(\alt,z)-5(h'(\alt,z))^2}{32h^{5/2}(\alt,z)}, \quad A_0(z)=1,\quad  A_j(z)\ll1.
\end{equation}
Applying  \cite[2.225.1, 2.224.5]{GR} we find that
\begin{multline}\label{xi eq}
\xi=-2\log(1+\sqrt{1-(1-\alt^2)z})+2\alt\log(\alt+\sqrt{1-(1-\alt^2)z})+\\+\log(z)-\alt\log(1-z)+
(1-\alt)\log(1-\alt^2).
\end{multline}
In order to determine $C_1$ and $C_2$ we compare two sides of equation \eqref{Y to W1W2} as $z\to0$. By \eqref{Y def} one has $Y(z)h^{1/4}(\alt,z)\sim z^{ir}$. As  a result, using  \eqref{xi eq} we infer that $C_2=0$ and
\begin{equation}\label{C1 def}
C_1=e^{ir\left(2\log2-2\alt\log(1+\alt)-(1-\alt)\log(1-\alt^2)\right)}.
\end{equation}
Finally,
\begin{multline}\label{2F1 LG asympt}
\HyG\left(1/4+ir(1-\alt),3/4+ir(1-\alt),1+2ir;z\right)=\\=
\frac{e^{2ir\Lo(\alt,z)}}{(1-(1-\alt^2)z)^{1/4}}
\left(1+\sum_{j=1}^{n-1}\frac{A_{j}(z)}{(ir)^j}\right)+O(r^{-n}),
\end{multline}
where $\Lo(\alt,z)$ is given by \eqref{2F1 main asymptI Lo func}, showing that \eqref{2F1 main asymptI} is proved for $z<1-\delta$.

The case when $z\to1$ is more complicated since in this situation the turning point coalesce with a pole as $\alpha\to0$. Asymptotic solutions of such differential equations have been investigated by Boyd and Dunster \cite{BD} and Dunster \cite{Dun90}. In particular, the differential equation \eqref{Y difeq} satisfies the conditions in \cite{Dun90}.  Nevertheless,  it is not obvious how to apply the results of \cite{Dun90} in order to study \eqref{Y def}. The problem is that in \cite{Dun90} only real solutions of the differential equation have been considered and our function $Y(z)$ is not real. Note that Dunster \cite{Dun90} applied his result to the study of Ferrers functions, which are not real in general. For this reason,  instead of individual functions Dunster studied their tricky combination \cite[(5.15)]{Dun90}, which ultimately is a real function.

Being unable to overcome this difficulty in our case, we decided to use  the saddle point method for proving \eqref{2F1 main asymptI} for values of $z$ that are close to one. Let us start with the case $z<1-\delta$.
Applying \cite[15.6.1]{HMF} one has
\begin{multline}
\HyG\left(1/4+ir(1-\alt),3/4+ir(1-\alt),1+2ir;z\right)=\\=
\frac{\Gamma(1+2ir)}{\Gamma\left(3/4+ir(1-\alt)\right)\Gamma\left(1/4+ir(1+\alt)\right)}\int_0^1
\frac{y^{-1/4}\exp(irf(\alt,y))}{(1-y)^{3/4}(1-zy)^{1/4}}dy,
\end{multline}
where
\begin{equation}
f(\alt,y)=(1-\alt)\log y+(1+\alt)\log(1-y)-(1-\alt)\log(1-zy).
\end{equation}
The saddle points (solutions of $\frac{\partial}{\partial y}f(\alt,y)=0$) are given by
$$y_{\pm}=\frac{1\pm\sqrt{1-(1-\alt^2)z}}{(1+\alt)z}.$$
For $z\to1$ and $\alt\to0$ both $y_{+}$ and $y_{-}$ tend to the end point of the segment of integration and coalesce if $z=1,\alt=0$. Such location of saddle points makes it impossible to apply the standard saddle point method.

For $0<z<1-\delta$ the point $y_{+}$ lies far away from the interval of integration, while the point $y_{-}$ is located inside the interval of integration and is not close to the end point.  In this case, the standard saddle point method yields the following asymptotic expansion:
\begin{multline}\label{2f1 sadpoint aprox}
\frac{\Gamma\left(3/4+ir(1-\alt)\right)\Gamma\left(1/4+ir(1+\alt)\right)}{\Gamma(1+2ir)}
\HyG\left(1/4+ir(1-\alt),3/4+ir(1-\alt),1+2ir;z\right)=\\=
\frac{y_{-}^{-1/4}e^{irf(\alt,y_{-})+\pi i/4}
\sqrt{2\pi}}{(1-y_{-})^{3/4}(1-zy_{-})^{1/4}\sqrt{rf''(\alt,y_{-})}}
\left(1+\sum_{j=1}^{N-1}\frac{c_j}{r^j}\right)+O(r^{-N-1/2}),
\end{multline}
where $c_j\ll 1.$ Direct computations show that
\begin{multline}\label{f at y-}
f(\alt,y_{-})=(1-\alt)\log(1-\alt^2)-2\log(1+\sqrt{1-(1-\alt^2)z})+2\alt\log(\alt+\sqrt{1-(1-\alt^2)z}),
\end{multline}
\begin{equation}\label{g at y-}
\frac{y_{-}^{-1/4}(1-y_{-})^{-3/4}}{(1-zy_{-})^{1/4}}=\frac{(1+\alt)^{1/4}(1+\sqrt{1-(1-\alt^2)z})}{(1-\alt)^{1/4}(\alt+\sqrt{1-(1-\alt^2)z})},
\end{equation}
\begin{equation}\label{f'' at y-}
f''(\alt,y_{-})=-2\frac{(1+\alt)\left(1+\sqrt{1-(1-\alt^2)z}\right)^2\sqrt{1-(1-\alt^2)z}}{(1-\alt)\left(\alt+\sqrt{1-(1-\alt^2)z}\right)^2}.
\end{equation}
Finally, using the Stirling formula we prove that
\begin{multline}\label{Stirling2}
\Gamma(\sigma+it)=\sqrt{2\pi}|t|^{\sigma-1/2}\exp(-\pi|t|/2)
\exp\left(i\left(t\log|t|-t+\frac{\pi t(\sigma-1/2)}{2|t|}\right)\right)\\\times
\left(1+\sum_{j=1}^{N-1}a_j/t^j+O(|t|^{-N})\right),
\end{multline}
which holds for $|t|\rightarrow\infty$ and a  fixed $\sigma$. Note that
\begin{multline}\label{2f1 sadpoint Gamma}
\frac{\Gamma(1+2ir)}{\Gamma\left(3/4+ir(1-\alt)\right)\Gamma\left(1/4+ir(1+\alt)\right)}=O(r^{-N+1/2})+\\+
\frac{(1+\alt)^{1/4}\sqrt{r}}{(1-\alt)^{1/4}\sqrt{\pi}}
e^{\pi i/4+ir\left(2\log2-(1+\alt)\log(1+\alt)-(1-\alt)\log(1-\alt)\right)}
\left(1+\sum_{j=1}^{N-1}\frac{c_j}{r^j}\right).
\end{multline}
Substituting \eqref{f at y-}, \eqref{2f1 sadpoint Gamma} to \eqref{2f1 sadpoint aprox} we complete the proof of \eqref{2F1 main asymptI}.
\section{The second case: $1-\delta<z<1$}
It is left to prove \eqref{2F1 main asymptI} for $1-\delta<z<1$. First, we transform \eqref{2f1(1/4 3/4 1)} in the same way as in
\cite[Lemma 4.6]{BF2mom}. Namely, using \cite[Eq. 16, p. 112]{BE} and \cite[Eq. 15.8.1]{HMF}, we show that
\begin{multline}\label{2f1 double transform}
\HyG\left(1/4+ir(1-\alt),3/4+ir(1-\alt),1+2ir;z\right)\\=
(1-z)^{-1/4-ir(1-\alt)}\HyG\left( \frac{1}{2}+2ir(1-\alt),\frac{1}{2}+2ir(1+\alt),1+2ir; \frac{\sqrt{1-z}-1}{2\sqrt{1-z}}\right)\\
=(1-z)^{-1/4-ir(1-\alt)}\left(\frac{\sqrt{1-z}+1}{2\sqrt{1-z}}\right)^{-2ir}
\HyG\left( \frac{1}{2}+2ir\alt,\frac{1}{2}-2ir\alt,1+2ir; \frac{\sqrt{1-z}-1}{2\sqrt{1-z}}\right).
\end{multline}
Let
\begin{equation}\label{Y def}
Y=\frac{1-\sqrt{1-z}}{2\sqrt{1-z}}
\end{equation}
Note that $Y\sim(1-z)^{-1/2}\to+\infty$ as $z\to1$. Applying \cite[15.6.1]{HMF} we infer
\begin{multline}
\HyG\left( \frac{1}{2}+2ir\alt,\frac{1}{2}-2ir\alt,1+2ir; -Y\right)=\\=
\frac{\Gamma(1+2ir)}{\Gamma\left(1/2+2ir\alt\right)\Gamma\left(1/2+2ir(1-\alt)\right)}\int_0^1
\frac{y^{-1/2+2ir\alt}(1-y)^{-1/2+2ir(1-\alt)}}{(1+Yy)^{1/2-2ir\alt}}dy.
\end{multline}
Next, we make the change of variables: $y:=1-e^{-x}$. As a result,
\begin{multline}\label{2f1 to Temme int}
\HyG\left(1/4+ir(1-\alt),3/4+ir(1-\alt),1+2ir;z\right)=\\=
\frac{(1-z)^{-1/4-ir(1-\alt)}(1+Y)^{-2ir}\Gamma(1+2ir)}{\Gamma\left(1/2+2ir\alt\right)\Gamma\left(1/2+2ir(1-\alt)\right)}\int_0^{\infty}
e^{2ir(-x+\alt\log q(x))}\frac{dx}{\sqrt{q(x)}},
\end{multline}
where
\begin{equation}\label{q(x) def}
q(x)=\left(1-e^{-x}\right)\left((Y+1)e^x-Y\right).
\end{equation}

In order  to obtain the asymptotic expansion for the integral in \eqref{2f1 to Temme int} we apply the saddle point method in the form developed by McKee, Sun, and Ye \cite{McKeeSunYe:2017}.
The saddle point of the integral \eqref{2f1 to Temme int} (the solution of $q(x)=\alt q'(x)$) is given by
\begin{equation}\label{Temme x0 def}
x_0=\log\frac{1+\sqrt{1-(1-\alt^2)z}}{(1-\alt)\left(1+\sqrt{1-z}\right)}.
\end{equation}
First, using the Taylor expansion, one can show that
 $\alt\ll x_0\ll\alt.$ Indeed, if $\alt^2 z / (1 - z)^{1 / 2} \ll 1$, then
\begin{multline}\label{x0 Taylor1}
x_0=
\log\frac{1+\sqrt{1-z}\sqrt{1+\alt^2z/(1-z)}}{1+\sqrt{1-z}}-\log(1-\alt)=\\=
\log\left(1+\frac{\alt^2z/2}{1-z+\sqrt{1-z}}+O\left(\frac{\alt^4z^2}{(1-z)^{3/2}}\right)\right)-\log(1-\alt)=\\=
\alt+\frac{\alt^2z/2}{1-z+\sqrt{1-z}}+O(\alt^2)+O\left(\frac{\alt^4z^2}{(1-z)^{3/2}}\right).
\end{multline}
Therefore, if $\alt \ll (1 - z)^{1 / 2}$, then $\alt \ll x_0 \ll \alt$. If $(1 - z)^{1 / 2} \ll \alt$, one has
\begin{multline}\label{x0 Taylor2}
x_0=\sqrt{1-z+\alt^2z}-\sqrt{1-z}+\alt+O\left(\alt^2\right)=
\alt+\frac{\alt^2z/\sqrt{1-z}}{1+\sqrt{1+\alt^2z/(1-z)}}+O\left(\alt^2\right),
\end{multline}
and again $\alt\ll x_0\ll\alt.$ Hence, there exist some absolute constants $D_0$ and $D_1$ such that $2\alt D_0<x_0<\alt D_1.$ To localize the contribution of the saddle point to the integral \eqref{2f1 to Temme int}, we split it into three parts using a partition of unity, $1=\chi_0(x)+\chi_{1}(x)+\chi_2(x)+\chi_3(x)$, consisting of the following smooth characteristic functions of intervals:
\begin{equation}\label{chi0 def}
\chi_0(x)=1\quad\hbox{for}\quad0<x<\alt D_0,\quad
\chi_0(x)=0\quad\hbox{for}\quad x>2\alt D_0,
\end{equation}
\begin{equation}\label{chi1 def}
\chi_{1}(x)=1\quad\hbox{for}\quad 2\alt D_0<x_0<\alt D_1,
\chi_{1}(x)=0\quad\hbox{for}\quad x<\alt D_0,\,x>2\alt D_1,
\end{equation}
\begin{equation}\label{chi2 def}
\chi_{2}(x)=1\quad\hbox{for}\quad 2\alt D_1<x<D_2\log^2(r\alt),\quad
\chi_{2}(x)=0\quad\hbox{for}\quad x<\alt D_1,\,x>2D_2\log^2(r\alt),
\end{equation}
\begin{equation}\label{chi3 def}
\chi_{3}(x)=1\quad\hbox{for}\quad x>2D_2\log^2(r\alt),\quad
\chi_{3}(x)=0\quad\hbox{for}\quad x<D_2\log^2(r\alt).
\end{equation}
As a result,
\begin{equation}\label{Temme int decompose1}
\int_0^{\infty}e^{2ir(-x+\alt\log q(x))}\frac{dx}{\sqrt{q(x)}}=\sum_{i=0}^{3}I_i(\alt,Y),\quad
I_i(\alt,Y)=\int_0^{\infty}g_i(x)e^{2\pi i f(x)}dx,
\end{equation}
where
\begin{equation}\label{fgi def}
f(x)=\frac{r}{\pi}(-x+\alt\log q(x)),\quad
g_i=\frac{\chi_i(x)}{\sqrt{q(x)}}.
\end{equation}
Estimating $I_3(\alt,Y)$ trivially with the use of \eqref{Y def}, \eqref{q(x) def}, and \eqref{chi3 def}, we get 
\begin{equation}\label{I3 est}
I_3(\alt,Y)\ll(1-z)^{1/4}(r\alt)^{-A}
\end{equation}
for an arbitrary $A>1.$ To estimate $I_0(\alt,Y)$, we integrate $n$ times by parts, as in \cite[(2.1)]{McKeeSunYe:2017}, obtaining
\begin{equation}\label{I0 int by parts1}
I_0(\alt,Y)=e^{2\pi i f(x)}\sum_{i=1}^{n}H_i(x)\Biggl|_{0}^{2\alt D_0}-
\int_0^{\infty}H'_n(x)e^{2\pi i f(x)}dx,
\end{equation}
where
\begin{equation}\label{Hi def}
H_1(x)=\frac{g_0(x)}{2\pi if'(x)},\quad H_i(x)=\frac{-H'_{i-1}(x)}{2\pi if'(x)}.
\end{equation}
Note that the contribution of the endpoint $2\alt D_0$ is zero due to the properties of $\chi_0$, see \eqref{chi0 def}.
It follows from \cite[(2.2)]{McKeeSunYe:2017} that for $k>1$,
\begin{equation}\label{Hk est1}
H_k(x)\ll\frac{g_0^{(k-1)}}{(f')^k}+\sum_{j=1}^{k-1}\frac{1}{(f')^{k+j}}\sum_{t=0}^{k-1-j}g_0^{(k-1-j-t)}
\sum_{\substack{2\le n_1,\ldots,n_j\le t+2\\ n_1+\ldots n_j=t+2j}}f^{(n_1)}\cdots f^{(n_j)}.
\end{equation}
Differentiating \eqref{fgi def}, we find that
\begin{equation}\label{f n der1}
f'(x)=\frac{r(\alt q'(x)-q(x))}{\pi q(x)},\quad
f^{(n)}=\frac{\alt r}{\pi} \left(\log q(x)\right)^{(n)}.
\end{equation}
By Fa\`{a} di Bruno's formula
\begin{equation}\label{FaadiBruno1}
h(g(x))^{(n)}\ll
\sum_{\substack{m_1+2m_2+\ldots nm_n=n\\m_j\ge0}}h^{(m_1+\ldots m_n)}(g)\prod_{j=1}^n\left(g^{(j)}\right)^{m_j}.
\end{equation}
Applying \eqref{FaadiBruno1} to $\left(\log q(x)\right)^{(n)}$, we obtain
\begin{equation}\label{f n der2}
f^{(n)}(x)\ll
\alt r \sum_{\substack{m_1+2m_2+\ldots nm_n=n\\m_j\ge0}}\prod_{j=1}^n\left(\frac{q^{(j)}(x)}{q(x)}\right)^{m_j}.
\end{equation}
It follows from \eqref{q(x) def} that 
\begin{equation}\label{q(x) def2}
q(x)=(Y+1)e^{x}+Ye^{-x}-1-2Y,
\end{equation}
\begin{equation}\label{q(x) deriv}
q^{(2j-1)}(x)=(Y+1)e^{x}-Ye^{-x},\quad
q^{(2j)}(x)=(Y+1)e^{x}+Ye^{-x},
\end{equation}
where $j\ge1$. Therefore, for small $x$, one has
\begin{equation}\label{q(x) deriv2}
q(x)\sim x+x^2Y,\quad
q^{(2j-1)}(x)\sim 1+xY,\quad
q^{(2j)}(x)\sim Y,
\end{equation}
and thus
\begin{equation}\label{q(x) deriv3}
\frac{q^{(2j-1)}(x)}{q(x)}\sim\frac{1}{x},\quad
\frac{q^{(2j)}(x)}{q(x)}\sim \frac{Y}{x(1+xY)}.
\end{equation}
Substituting \eqref{q(x) deriv3} into \eqref{f n der2}, we find that for $n\ge2$
\begin{equation}\label{f n der3}
f^{(n)}(x)\ll
\alt r \sum_{\substack{m_1+2m_2=n\\m_j\ge0}}
\frac{Y^{m_2}}{x^{m_1+m_2}(1+xY)^{m_2}}\ll\frac{\alt r }{x^n}.
\end{equation}
The last inequality is trivial in the case when $xY>1$, and if $xY<1$,
\begin{equation*}
\sum_{\substack{m_1+2m_2=n\\m_j\ge0}}
\frac{Y^{m_2}}{x^{m_1+m_2}(1+xY)^{m_2}}\ll\sum_{\substack{m_1+2m_2=n\\m_j\ge0}}\frac{Y^{m_2}}{x^{n-m_2}}\ll\frac{1}{x^n}.
\end{equation*}
Also, for $0<x<2\alt D_0$, combining \eqref{f n der1} and \eqref{q(x) deriv3}, we obtain
\begin{equation}\label{f n der4}
f'(x)\sim r\left(\frac{\alt}{x}-1\right)\sim\frac{\alt r}{x}.
\end{equation}
Next, to be able to estimate \eqref{Hk est1}, we must investigate $g_0^{(k)}(x)$. It follows from \eqref{fgi def} and \eqref{chi0 def} that
\begin{equation}\label{g0 der1}
g_0^{(k)}(x)\ll \chi_0(x)(q(x)^{-1/2})^{(k)}+\sum_{1\le n\le k}\left(\frac{1}{\sqrt{q(x)}}\right)^{(n)}\frac{1}{(\alt D_0)^{k-n}}.
\end{equation}
Arguing in the same way as in the proof of \eqref{f n der3}, we obtain
\begin{equation}\label{g0 der2}
\left(\frac{1}{\sqrt{q(x)}}\right)^{(n)}\ll
\frac{1}{x^n\sqrt{x+x^2Y}},
\end{equation}
and thus 
\begin{equation}\label{g0 der3}
g_0^{(k)}(x)\ll \frac{1}{x^k\sqrt{x+x^2Y}}.
\end{equation}
With the help of \eqref{f n der3}, \eqref{f n der4}, and \eqref{g0 der3}, we estimate \eqref{Hk est1} as
\begin{equation}\label{Hk est2}
H_k(x)\ll
\frac{x}{(\alt r)^k\sqrt{x+x^2Y}}+
\sum_{j=1}^{k-1}\frac{x^{k+j}}{(\alt r)^{k+j}}\sum_{t=0}^{k-1-j}\frac{1}{x^{k-1-j-t}\sqrt{x+x^2Y}}
\frac{(\alt r)^j}{x^{t+2j}}\ll \frac{x}{(\alt r)^k\sqrt{x+x^2Y}}.
\end{equation}
Therefore, the contribution of the endpoint $x=0$ to \eqref{I0 int by parts1} is zero and
\begin{equation}\label{I0 int by parts2}
I_0(\alt,Y)\ll \sup|H_{n+1}(x)|\ll \sup\frac{\sqrt{x}}{(\alt r)^{n+1}\sqrt{1+xY}}\ll\frac{1}{(\alt r)^{n+1}\sqrt{Y}}.
\end{equation}
Arguing in the same way, one can obtain the same estimate for $I_2(\alt,Y)$. Using this result, \eqref{I0 int by parts2}, and \eqref{I3 est}, we obtain from \eqref{Temme int decompose1} that
\begin{equation}\label{Temme int decompose2}
\int_0^{\infty}e^{2ir(-x+\alt\log q(x))}\frac{dx}{\sqrt{q(x)}}=I_1(\alt,Y)+O\left(\frac{1}{(\alt r)^{n+1}\sqrt{Y}}\right),
\end{equation}
for any  $n\ge1$. 
It remains to obtain an asymptotic formula for $I_1(\alt,Y)$. To do this, we apply \cite[Theorem 1.2]{McKeeSunYe:2017}, which requires determining $N,M,T,U$ such that for $s=0,\ldots,2n+1$ and $r=2,\ldots,2n+3$
\begin{equation}\label{g1 and f deriv est1}
|g_{1}^{(s)}(x)|\ll\frac{U}{N^{s}}, \quad
|f^{(r)}(x)|\ll\frac{T}{M^r} \quad\hbox{and}\quad
|f^{(2)}(x)|\gg\frac{T}{M^2}.
\end{equation}
Similarly to \eqref{g0 der3}, we obtain
\begin{equation}\label{g1 der}
g_1^{(k)}(x)\ll \frac{1}{\alt^k\sqrt{\alt+\alt^2Y}},
\end{equation}
and thus $U=(\alt+\alt^2Y)^{-1/2}$ and $N=\alt.$  
According to \eqref{fgi def} and \eqref{q(x) deriv2}
\begin{equation}\label{f 2der1}
f''(x)=\frac{r\alt}{\pi}\frac{q''(x)q(x)-(q'(x))^2}{q^2(x)}\gg r\alt\frac{q''(x)q(x)-(q'(x))^2}{\alt^2(1+\alt Y)^2}.
\end{equation}
Applying \eqref{q(x) def2} and  \eqref{q(x) deriv}, we find that
\begin{multline}\label{f 2der2}
q''(x)q(x)-(q'(x))^2=4Y(Y+1)-(Y+1)(2Y+1)e^x-Y(2Y+1)e^{-x}\\=
-1-(2Y+1)^2\sum_{\substack{n=1\\n\equiv0\Mod{2}}}^{\infty}\frac{x^n}{n!}-
(2Y+1)\sum_{\substack{n=1\\n\equiv1\Mod{2}}}^{\infty}\frac{x^n}{n!}\\=-1-(2Y+1)x-(2Y+1)^2x^2/2+O(Y^2x^3+Yx^2).
\end{multline}
Substituting \eqref{f 2der2} into \eqref{f 2der1}, we show that
\begin{equation}\label{f 2der3}
f''(x)\gg r\alt\frac{1+\alt Y+\alt^2Y^2}{\alt^2(1+\alt Y)^2}\gg\frac{r\alt}{\alt^2}.
\end{equation}
It follows from \eqref{f n der3} and \eqref{f 2der3}  that \eqref{g1 and f deriv est1} holds with $M=\alt$ and $T=r\alt$.
Thus, using  \cite[Theorem 1.2]{McKeeSunYe:2017} to  evaluate the integral $I_1(\alt,Y)$, we obtain
\begin{multline}\label{Temme int decompose33}
\int_0^{\infty}e^{2ir(-x+\alt\log q(x))}\frac{dx}{\sqrt{q(x)}}=
\frac{\sqrt{\pi}e^{-\pi i/4}e^{2ir(-x_0+\alt\log q(x_0))}}{\sqrt{r\alt}\sqrt{(q'(x_0))^2/q(x_0)-q''(x_0)}}
\left(1+\sqrt{q(x_0)}\sum_{j=1}^{n}d_j(\alt,r,Y)\right)
+\\
+O\left(\frac{1}{(\alt r)^{n+1}\sqrt{Y}}+\frac{\alt}{(r\alt)^{n+1}\sqrt{\alt+\alt^2Y}}\right),
\end{multline}
where (see \cite[(1.9)-(1.11)]{McKeeSunYe:2017})
\begin{equation}\label{sp coef1}
d_k(\alt,r,Y)\ll\frac{g_1^{(2k)}}{(f'')^k}+\sum_{l=0}^{2k-1}\frac{g_1^{(l)}}{(f'')^k}\sum_{j=1}^{2k-l}\frac{1}{(f'')^j}
\sum_{\substack{n_1+m_2+\ldots n_j=2k-l+2j\\3\le n_1+\ldots +n_j\le2n+3}}f^{(n_1)}\ldots f^{(n_j)},
\end{equation}
with all derivatives being evaluated at the point $x_0$. Using \eqref{f 2der3}, \eqref{f n der3} and \eqref{g1 der}, we infer
\begin{equation}\label{sp coef2}
d_k(\alt,r,Y)\ll\frac{1}{(r\alt)^k\sqrt{\alt+\alt^2Y}}.
\end{equation}
From \eqref{sp coef2}, applying \eqref{q(x) deriv2}, we obtain that
\begin{equation}\label{sp coef3}
\sqrt{q(x_0)}d_k(\alt,r,Y)\ll\frac{1}{(r\alt)^k}.
\end{equation}
Using \eqref{sp coef3}, we rewrite \eqref{Temme int decompose33} as
\begin{multline}\label{Temme int decompose3}
\int_0^{\infty}e^{2ir(-x+\alt\log q(x))}\frac{dx}{\sqrt{q(x)}}=
\frac{\sqrt{\pi}e^{-\pi i/4}e^{2ir(-x_0+\alt\log q(x_0))}}{\sqrt{r\alt}\sqrt{(q'(x_0))^2/q(x_0)-q''(x_0)}}
\left(1+\sum_{j=1}^{n}\frac{c_j(\alt,r,z)}{(r\alt)^j}\right)
+\\+O\left(\frac{1}{(\alt r)^{n+1}\sqrt{Y}}\right),
\end{multline}
with $c_j(\alt,r,z)\ll1.$


Straightforward computations show that
\begin{equation}\label{Temme mainterm2}
\frac{(q'(x_0))^2}{q(x_0)}-q''(x_0)=
\frac{q(x_0)}{\alt^2}-q''(x_0)=\frac{\sqrt{1-(1-\alt^2)z}}{\alt\sqrt{1-z}},
\end{equation}
\begin{equation}\label{Temme mainterm3}
-x_0+\alt\log q(x_0)=-\log\frac{1+\sqrt{1-(1-\alt^2)z}}{(1-\alt)(1+\sqrt{1-z})}+
\alt\log\frac{\alt+\sqrt{1-(1-\alt^2)z}}{(1-\alt^2)\sqrt{1-z}}+\alt\log\alt.
\end{equation}
Substituting \eqref{Temme mainterm2} and \eqref{Temme mainterm3} to \eqref{Temme int decompose3} we conclude that the main term in \eqref{Temme int decompose3} is given by
\begin{equation}\label{2f1 to Temme int3}
\frac{e^{-\pi i/4}\sqrt{\pi}}{\sqrt{r}}\frac{(1-z)^{1/4}}{(1-(1-\alt^2)z)^{1/4}}e^{2ir\left(
\alt\log\alt-\log\frac{1+\sqrt{1-(1-\alt^2)z}}{(1-\alt)(1+\sqrt{1-z})}+
\alt\log\frac{\alt+\sqrt{1-(1-\alt^2)z}}{(1-\alt^2)\sqrt{1-z}}
\right)}.
\end{equation}
Substituting \eqref{Temme int decompose3} and \eqref{2f1 to Temme int3} to  \eqref{2f1 to Temme int}, applying
\begin{equation}
\frac{\Gamma(1+2ir)}{\Gamma\left(1/2+2ir\alt\right)\Gamma\left(1/2+2ir(1-\alt)\right)}
=\frac{\sqrt{r}}{\sqrt{\pi}}e^{\pi i/4+2ir(-\alt\log\alt-(1-\alt)\log(1-\alt))}\left(1+\sum_{j=1}^{N-1}\frac{d_j}{(r\alt)^j}\right)+O\left(\frac{\sqrt{r}}{(r\alt)^{N}}\right),
\end{equation}
where $d_j\ll1$ and using the fact that $Y=\frac{1-\sqrt{1-z}}{2\sqrt{1-z}}$, we obtain
\begin{multline}\label{2f1 Temme asympt1}
\HyG\left(1/4+ir(1-\alt),3/4+ir(1-\alt),1+2ir;z\right)=\\=
\frac{e^{2irF(\alt,z)}}{(1-(1-\alt^2)z)^{1/4}}\left(1+\sum_{j=1}^{N-1}\frac{c_j(\alt,r,z)}{(r\alt)^j}+O\left(\frac{\sqrt{r}}{(r\alt)^{N}}\right)\right),
\end{multline}
where $c_j(\alt,r,z)\ll1$ and
\begin{multline}\label{2f1 to Temme mt F def}
F(\alt,z)=-\log\frac{1+\sqrt{1-(1-\alt^2)z}}{(1-\alt)(1+\sqrt{1-z})}+
\alt\log\frac{\alt+\sqrt{1-(1-\alt^2)z}}{(1-\alt^2)\sqrt{1-z}}-\\-
\log\frac{1+\sqrt{1-z}}{2\sqrt{1-z}}-(1-\alt)\log\sqrt{1-z}-(1-\alt)\log(1-\alt).
\end{multline}
One can easily check that $F(\alt,z)$ coincides with $\Lo(\alt,z)$  defined in \eqref{2F1 main asymptI Lo func}. This completes the proof of \eqref{2F1 main asymptI}.


\section{Numerical examples}

In this section, we provide some numerical examples related to Theorem \ref{prop: 2F1 asympt1}.
More precisely, we compute for some values of $r,$ $\alpha$ and $z$ the left-hand side of \eqref{2F1 main asymptI}, given by
$$F(r,\alpha,z):=\HyG\left(1/4+ir(1-\alt),3/4+ir(1-\alt),1+2ir;z\right),$$ and compare it
with the first summand on the right-hand side of \eqref{2F1 main asymptI}:
$$R(r,\alpha,z):=\frac{\exp\left(2ir\Lo(\alt,z)\right)}{(1-(1-\alt^2)z)^{1/4}}.$$

The results are given in the Tables below.

\begin{center}
 \begin{tabular}{ |c|c|c| }
\hline
Function & Approximation & $|$ Relative error $|$ \\ \hline
 &  &  \\
$F(r,\alpha,z)$  & $2.611880247+0.5174226366i$ &  \\
&  &  \\
$R(r,\alpha,z)$ & $2.611802650+0.5170096337i$ & $0.000157843$  \\
&  &  \\ \hline
\end{tabular}
\captionof{table}{The case $r=100$, $\alpha=0.1$, $y=0.99$.}\label{case4}
\end{center}


\begin{center}
 \begin{tabular}{ |c|c|c| }
\hline
Function & Approximation & $|$ Relative error $|$ \\ \hline
 &  &  \\
$F(r,\alpha,z)$  & $-2.595771772-1.792471289i$ &  \\
&  &  \\
$R(r,\alpha,z)$ & $-2.587392554-1.804512206i$ & $0.004650302$  \\
&  &  \\ \hline
\end{tabular}
\captionof{table}{The case $r=100$, $\alpha=0.1$, $y=0.9999$.}\label{case4}
\end{center}


\begin{center}
 \begin{tabular}{ |c|c|c| }
\hline
Function & Approximation & $|$ Relative error $|$ \\ \hline
 &  &  \\
$F(r,\alpha,z)$  & $1.000002393+0.004050009913i$ &  \\
&  &  \\
$R(r,\alpha,z)$ & $1.000002393+0.004050054037i$ & $7.576766037*10^{-9}$  \\
&  &  \\ \hline
\end{tabular}
\captionof{table}{The case $r=100$, $\alpha=0.1$, $y=0.00001$.}\label{case4}
\end{center}


\begin{center}
 \begin{tabular}{ |c|c|c| }
\hline
Function & Approximation & $|$ Relative error $|$ \\ \hline
 &  &  \\
$F(r,\alpha,z)$  & $-2.258587650-2.168919043i$ &  \\
&  &  \\
$R(r,\alpha,z)$ & $-2.269908138-2.1575933584i$ & $0.005113882$  \\
&  &  \\ \hline
\end{tabular}
\captionof{table}{The case $r=100$, $\alpha=0.02$, $y=0.99$.}\label{case4}
\end{center}


\begin{center}
 \begin{tabular}{ |c|c|c| }
\hline
Function & Approximation & $|$ Relative error $|$ \\ \hline
 &  &  \\
$F(r,\alpha,z)$  & $-3.328488374+5.815264147i$ &  \\
&  &  \\
$R(r,\alpha,z)$ & $-3.380721684+5.770084005i$ & $0.010307062$  \\
&  &  \\ \hline
\end{tabular}
\captionof{table}{The case $r=100$, $\alpha=0.02$, $y=0.9999$.}\label{case4}
\end{center}


\end{document}